\documentclass[12pt, oneside]{amsart}
\usepackage{a4wide,amsmath,amssymb,amscd,verbatim}
\input xypic

\newtheorem{prop}[equation]{Proposition}

\newtheorem{thm}[equation]{Theorem}
\newtheorem{cor}[equation]{Corollary}
\newtheorem{lem}[equation]{Lemma}

\theoremstyle{definition}

\newtheorem{defn}[equation]{Definition}
\newtheorem{rem}[equation]{Remark}
\numberwithin{equation}{section}

\newcommand{\letbe}{\!\stackrel{\text{def}}{=}}

\newcommand{\col}{\operatorname{colim}}

\newcommand{\hoc}{\operatorname{hocolim}}

\newcommand{\bC}{\mathbb{C}}
\newcommand{\bF}{\mathbb{F}}

\newcommand{\bQ}{\mathbb{Q}}
\newcommand{\bR}{\mathbb{R}}
\newcommand{\bZ}{\mathbb{Z}}

\newcommand{\cat}[1]{\mbox{\sc #1}}

\newcommand\Z{\bZ}
\newcommand\z{\Z}
\newcommand\R{\bR}
\newcommand\q{\bQ}
\newcommand\Q{\bQ}
\newcommand \fp{\mathbb F_p}
\newcommand \lra{\longrightarrow}
\newcommand \p{^\wedge_p}

\newcommand \ra{\rightarrow}
\def \larrow#1{\,\stackrel{#1}\lra\,}

\def \CP{\mathcal{P}}

\newcommand{\const}{\operatorname{const}}

\newcommand{\catk}{\cat{cat}(K)}
\newcommand{\catkop}{\catk^{op}}
\newcommand{\catcat}{\cat{cat}}
\newcommand{\ab}{\cat{ab}}
\newcommand{\Ab}{\cat{ab}}

\newcommand{\nz}{\newline}

\newcommand{\catkx}{\catcat(K^\times)}
\newcommand{\catc}{\cat{c}}
\newcommand{\catd}{\cat{d}}
\newcommand{\op}{^{op}}
\newcommand{\catkxop}{\catkx^{op}}
\newcommand{\Top}{\cat{Top}}
\newcommand{\Topx}{\cat{Top}^+}
\newcommand{\da}{\downarrow}
\newcommand{\Dim}{\operatorname{dim}}%

\newcommand{\soc}{\operatorname{soc}}%
\newcommand{\link}{\operatorname{link}}

\newcommand{\map}{\operatorname{map}}
\newcommand{\Hom}{\operatorname{Hom}}
\newcommand{\Tor}{\operatorname{Tor}}

\newcommand{\st}{\operatorname{st}}
\newcommand{\cst}{\operatorname{cst}}
\newcommand{\bst}{\operatorname{F}}
\newcommand{\Tot}{\operatorname{Tot}}

\newcommand{\id}{\operatorname{id}}
\begin{document}
\bibliographystyle{plain}
%\centerline{\hfill \lebn}
%\vspace*{1cm}

\title{Cohen-Macaulay and Gorenstein complexes from a topological 
point of view
}

\author{Dietrich Notbohm}
\address{Department of Mathematics, University of
Leicester, University Road, Leicester LE1~7RH, England}
\email{dn8@mcs.le.ac.uk}


\begin{abstract}
The main invariant to study the combinatorics of a 
simplicial complex $K$ is the associated face ring or Stanley-Reisner algebra.
Reisner respectively Stanley explained in which sense 
Cohen-Macaulay and Gorenstein
properties of the face ring are reflected by geometric and/or 
combinatoric properties of the simplicial complex. We give a new proof for 
these result by homotopy theoretic methods and constructions. 
Our approach is based on 
 ideas used very successfully in the analysis of the homotopy theory
of classifying spaces.
\end{abstract}

\maketitle

\section{Introduction} \label{introduction}

The main tool and invariant for understanding the combinatorics 
of a finite  simplicial complex is the associated face ring or
Stanley-Reisner algebra which is a quotient of a polynomial algebra 
generated by the vertices.
It is of interest to which extend
algebraic properties of the face ring are reflected by
combinatorial or geometric
properties of the simplicial complex. 
For example, Reisner
characterized all simplicial complexes whose face ring is Cohen-Macaulay
\cite{re}. And 
Stanley proved a similar result for Gorenstein face rings \cite{st}.
In this paper we will look at these result with the eyes of a topologist
and reprove both results with methods and ideas from homotopy theory,
in particular from the homotopy theory of classifying spaces.
For the topological proof we introduce and discuss some 
new spaces associated with 
simplicial complexes,
which, as we feel, deserve interest in their own right.

Let  $K$ be an abstract simplicial complex with m vertices given by the set 
$V=\{1,...,m\}$
That is, $K=\{\sigma_1, ..., \sigma_r\}$ consists of a finite set of faces 
$\sigma_i \subset V$, which  is closed
with respect to  formation of subsets.
The dimension
of $K$ is denoted by $\Dim K = n-1$. That is every
face $\sigma$ of $K$ has
order $\sharp\sigma\leq n$ and there exists a face $\mu$ of
order $\sharp\mu=n$. We consider the empty set $\emptyset$ as a face of $K$.

For a commutative ring $R$ with unit we denote by $R(K)$ the associated 
Stanley-Reisner algebra
of $K$ over the ring $R$. It is the quotient
$R[V]/(v_\sigma : \sigma \not\in K)$, where
$R[V]\letbe R[v_1,...v_m]$ is the polynomial algebra on $m$-generators
and $v_\sigma\letbe \prod_{j\in \sigma} v_j$. 
We can think of $R(K)$ as a
graded object. Since we want to bring topology into the game we will choose 
the topological grading and  give 
the generators of $R(K)$ and $R[V]$ the degree 2. 

Each abstract simplicial complex $K$ has a geometric realization, 
denoted by $|K|$. Let $\Delta(V)$ denote the full simplicial complex
whose faces are given by all subsets of $V$. One realization
of $\Delta(V)$ is given by the standard $m-1$-dimensional simplex 
$|\Delta(V)|\letbe \{\Sigma_{i=1}^m t_ie_i : 0\leq t_i\leq 1, \Sigma_i t_i=1\}
\subset \R^m$, where $e_1,...,e_m$ denotes the standard basis.
And a topological realization of $K$ is given by the subset 
$|K|\subset |\Delta(V)|$ defined by the subset relation $K\subset \Delta(V)$.
We define the homology $H_*(K)$ and cohomology $H^*(K)$ of $K$ as 
the homology $H^*(|K|)$ and cohomology $H^*(|K|)$ of the 
topological realization.

Before we can state the theorems of Reisner and Stanley, we have 
to recall some notions.
We call a simplicial complex $K$ Cohen-Macaulay or Gorenstein over a 
field $\bF$ if 
$\bF(K)$ is a Cohen-Macaulay or Gorenstein algebra over $\bF$.
We call $K$ a Gorenstein$^*$ complex if it is Gorenstein and if $V$ equals
the union of all minimal missing simplices of $K$. A subset $\mu \subset V$ 
is minimal missing if $\mu\not\in K$ and for each $\sigma \subsetneq \mu$
we have $\sigma \in K$. Moreover,
for any face $\sigma \in K$, the link of $\sigma$ is defined 
as the simplicial complex 
$$
\link_K(\sigma)\letbe \{\tau\setminus \sigma : \sigma \subset \tau \in K\}.
$$
Now the theorems of Reisner and Stanley read as follows:

\begin{thm} (Reisner \cite{re}) \label{cmresult}
Let $\bF$ be a field and $K$ a simplicial complex. Then, $K$ is 
Cohen-Macaulay over $\bF$
if and only if
for each face $\sigma$ of $K$, including the empty face,
$$
\widetilde H^i(\link_K(\sigma);\bF) = 0 \text{ for } i< \Dim \link_K(\sigma)
$$ 
\end{thm}

\begin{thm} (Stanley \cite{st}) \label{gorensteinresult}
Let $\bF$ be a field and $K$ a simplicial complex. Then, $K$ is 
Gorenstein$^*$ over $\bF$ if and only if 
$$
\widetilde H^i(\link_K(\sigma);\bF)\cong
\left\{
\array{ll}
\bF & \text{ for } i=\Dim \link_K(\sigma) \\
0 & \text{ for } i\neq \Dim \link_K(\sigma)
\endarray
\right.
$$
\end{thm}

For a definition of Cohen-Macaulay and Gorenstein properties see
either \cite{brhe}, \cite{st} or Section \ref{dcmg}.
There also exists a version of the second statement which deals 
with general Gorenstein complexes. But for simplification we
formulated the result for Gorenstein$^*$ complexes.

Reisner used methods of commutative algebra, 
in particular the machinery
of local cohomology of modules, to prove his theorem.
Another proof by Hochster (unpublished, see \cite{st})
is based on similar methods and a detailed analysis 
of the Poincar\'e series of $R(K)$. Similar ideas were used by Stanley
in his approach towards Gorenstein complexes.
In this paper we want to give a different proof for both results. Our proof
is based
on topological constructions related to and based on 
topological interpretations
of the
combinatorial data of $K$. For example, there exists a topological space $c(K)$
such that $H^*(X;\z) \cong \z(K)$ as a graded algebra. These spaces can be 
constructed as the Borel construction of toric manifolds 
\cite{daja}, as a (pointed) colimit of a particular diagram
\cite{bupa} or as the homotopy colimit of the same diagram
\cite{noraa}. This last construction is the most appropriate for doing 
homotopy theory and is the one which we will use in this paper. 
Using the homotopy colimit construction, if $\Dim K=n-1$, 
one can construct a very interesting non trivial map
$f:c(K) \lra BU(n)$ \cite{norac}. 
The Chern classes of the associated  vector bundle are given by the 
elementary symmetric polynomials in the generators of 
$\Z(K)$ (see Section \ref{vectorbundle}).
The homotopy fibre $X_K$ is a finite $CW$-complex, which  
contains a large amount of information about the associated 
Stanley-Reisner algebra. 

\begin{thm} \label{charfibre} $\ \ \ \ $
\nz
(i) A simplicial complex $K$ is Cohen-Macaulay 
over $\bF$,  if and only if $H^*(X_K;\bF)$
is concentrated in even degrees.
\nz
(ii) $K$ is Gorenstein over $\bF$ if and only if 
$H^*(X_K;\bF)$ is a Poincar\'e duality algebra and concentrated 
in even degrees. 
\end{thm}

This theorem translates Cohen-Macaulay and Gorenstein properties of $\bF(K)$ 
into conditions on $X_K$ and is the key result necessary for our 
proof of the results
of Reisner and Stanley.
 
The paper is organised as follows. In the next two sections we 
provide the basic 
topological ingredients necessary for the proof of
Theorem \ref{cmresult} and Theorem \ref{gorensteinresult}. In particular,
we recall the above mentioned homotopy colimit construction for the 
space $c(K)$ in Section \ref{hodec} and discuss the map $f_K : c(K) \lra BU(n)$
in Section \ref{vectorbundle}.
In Section \ref{dcmg} we provide 
definitions for Cohen-Macaulay and Gorenstein 
algebras appropriate for our purpose, express these properties in terms 
of the homotopy fibre $X_K$ and prove Theorem \ref{charfibre}. 
In Section \ref{hofipo} we discuss homotopy fixed point sets and study 
them for particular actions of elementary abelian 
groups on $X_K$. 
The final three sections are devoted to a proof
of Theorem \ref{cmresult} and Theorem \ref{gorensteinresult}.  

We will fix the following notation throughout.
$K$ always denotes a $(n-1)$-dimensional abstract finite simplicial complex
with $m$-vertices. The set of vertices will be denoted by 
$V\letbe V_K\letbe \{1,...,m\}$. We denote by $\ab$ the category of 
abelian groups, by
$\Top$ the category of topological spaces and by $\Topx$ the category 
of pointed topological spaces.  Mostly, we are working over  
commutative rings with unit  or  fields.
In particular, $R$ will always denote such a commutative ring and $\bF$ a field.
When we deal with torsion groups, we will use the topological convention 
for the grading. Projective resolutions are considered as non positively graded 
cochain complexes
and torsion groups are non positively graded objects denoted by 
$\Tor^j_A(M,N)$, where $j\leq 0$.

It is my pleasure to thank N. Ray and T. Panov for introducing me to this 
subject as well as for plenty of helpful discussions.

\bigskip

\section{Pointed and unpointed homotopy colimits} \label{hodec}

Given a category $\catc$ and functor $F : \catc \lra \Top$,
the colimit $\col_{\catc} F$ behaves particularly poorly in the context of
homotopy theory, since object-wise equivalent diagrams may
well have homotopy inequivalent colimits. The standard procedure for
dealing with this situation is to introduce the left derived
functor, known as the {\it homotopy colimit}. Following
\cite{hovo}, for example, $\hoc_{\catc} F$ may be described by the
two-sided bar construction.
In a similar fashion, we can construct a pointed homotopy colimit
$\hoc^+_{\catc} G$ for a functor $G:\catc \lra \Topx$. Composing $G$ 
with the forgetful functor $\phi: \Topx \lra \Top$ induces an identity
$\col_{\catc} \phi G = \col^+_{\catc} G$ and a cofibre sequence
$$
B\catc \lra \hoc_{\catc} \phi G \lra \hoc^+_{\catc} G
$$ 
\cite{boka}. Here, $B\catc$ is the classifying space of the
category $\catc$, that is the topological realization of the 
nerve $N(\catc)$. For details of the homotopy colimit construction 
see \cite{boka} \cite{drfa} or \cite{dw}.

The cohomology of (pointed) homotopy colimits can be calculated with the 
help of the Bousfield-Kan spectral sequence \cite{boka}. 
And this tool is more 
important for our purpose than the actual construction of the homotopy colimit.
In both cases, this is a first quadrant spectral sequence and has 
in the case of the unpointed homotopy colimit the form
$$
E_2^{i,j}\letbe  \lim{}^i_{\catc^{op}} H^j(F) 
\Longrightarrow H^{i+j}(\hoc_{\catc} F)
$$
and in the case of pointed homotopy colimit the form
$$
E_2^{i,j}\letbe  \lim{}^i_{\catc^{op}}\widetilde{H}^j(F) 
\Longrightarrow \widetilde{H}^{i+j}(\hoc_{\catc} F)
$$
In both cases, differentials $d_r : E_r^{*,*}\lra E_r^{*,*}$ 
have the degree $(r,1-r)$.

 Let $\phi : \catc^{op} \lra \ab$ be a
(covariant) functor, e.g. $\phi\letbe H^*(F)$.
Higher derived limits of $\phi$ can be thought of as the cohomology groups of a 
certain cochain complex 
$(C^*(\catc,F), \delta)$. The groups
are defined as
$$
C^r(\catc;F)\letbe \prod_{c_o\ra c_1 \ra... \ra c_r} F(c_0) \ \ \
\text{ for } r \geq 0 .
$$
Here, the morphism are morphisms in $\catc$ and not in $\catc^{op}$.
The differential $\delta : C^n(\catc;F) \lra C^{n+1}(\catc;F)$
is given by the alternating sum
$\sum_{k=0}^r (-1)^k \delta^k$ where $\delta_k$ is 
defined on $u\in C^n(\phi)$ by
$$
\delta^k(u)(c_0\ra\dots\ra c_{n+1})\letbe
\begin{cases}
u(c_0\ra\dots\ra\widehat{c}_k
\ra\dots\ra c_{n+1})&\text{for $k\neq 0$}\\
\phi(c_o\ra c_1)u(c_1\ra\dots
\ra c_n)&\text{for $k=0$}.
\end{cases}
$$
We may and will replace $C^*(\catc;F)$ by it's quotient
$N^*(\catc,F)$ of normalised cochains, given by the product over
chains
$c_0 \ra c_1 ... \ra c_r$ with distinct objects. 

In most cases we consider, the Bousfield-Kan 
spectral sequence  has a particular 
simple form. All higher limits will vanish. Following Dwyer \cite{dw}
we turn this property into the following definition.

\begin{defn}
Let $F : \catc \lra \Top$ be a functor and $R$  a commutative ring with unit. 
We call a map 
$f: \hoc_{\catc} F \lra Y$ a sharp $R$-homology decomposition,
if $\lim_{\catc^{op}}^i H^*(F;R)=0$ for $i\geq 1$ and if 
$f$ induces an isomorphism $\lim_{\catc^{op}} H^*(F;R)\cong H^*(Y;R)$.
If $F$ takes values in $\Top^+$, then we replace cohomology by reduced 
cohomology. 
\end{defn} 

Using the subset relation on  the faces, a simplicial
complex $K$ can be
interpreted as a poset and therefore as a
category which we denote by $\catk$. This category contains the empty set
$\emptyset$ as an initial object. If we want to exclude $\emptyset$ we
denote this by
$K^\times$ respectively by $\catkx$. In these cases the classifying space
$B\catk$ is equivalent to the cone of the topological realization $|K|$ and 
$B\catkx\simeq |K|$.

For a pointed topological space $Y$
we can define functors
$$
Y^K : \catk \lra \Topx, \ \ Y^K : \catk \lra \Top,
$$
which assigns the Cartesian product $Y^\sigma$
 to each face $\sigma$ of $K$. The value of $Y^K$ on $\sigma \subset \tau$
is the inclusion $Y^\sigma \subset Y^\tau$ where the superfluous 
coordinates are set to $*$. We note that $Y^K(\emptyset)\letbe *$. Moreover, 
this functor 
comes with a natural transformation $Y^K\lra 1_{Y^V}$ where $1_{Y^V}$ 
denotes the constant functor mapping each face to $Y^V=Y^m$ and 
each morphism to 
the identity map.

Since $\catk$ has an initial object, the classifying space $B\catk$ 
is contractible \cite{boka}  
and the above cofibre sequence degenerates to a 
homotopy equivalence
$$
\hoc_{\catk} Y^K\larrow{\simeq} \hoc^+_{\catk} Y^K.
$$

We want to specialise further. Let $T\letbe S^1$ denote the 1-dimensional torus
and $BT=\bC P^\infty$ the classifying 
space of $T$ respectively the infinite dimensional complex 
projective space. 
For the functor $BT^K :\catk \lra \Top$, in fact for all functors of 
the form $Y^K$,  the 
higher derived limits of 
the Bousfield-Kan spectral sequence vanish and the spectral sequence collapses
at the $E_2$-term. Only the actual inverse limit contributes 
something non trivial, namely the associated Stanley-Reisner algebra
\cite{noraa}. Defining $c(K)\letbe \hoc_{\catk} BT^K$ we can formulate this as
follows.

\begin{thm} \cite{noraa} \label{firstdecomposition}
\nz
(i) $H^*(\hoc_{\catk} BT^K; R)\cong R(K)$.
\nz
(ii)
$\hoc_{\catk} BT^K \lra c(K)$ is a sharp $R$-homology decomposition as 
well as 
\nz
$\hoc^+_{\catk} BT^K  \lra c(K)$.
\nz
(iii) The natural map $c(K) \lra BT^m$ realizes the algebra map 
$\Z[V] \lra \z(K)$.
\end{thm}

\begin{proof}
Part (i) and (iii) and the first half of Part (ii)are already 
proven in \cite{noraa}. 
The second claim of Part (ii) follows from 
the equivalence between the pointed and unpointed homotopy colimit, 
from the fact 
that reduced cohomology is a natural retract from cohomology 
and from part (i).
\end{proof}

For our purpose we will also need another homology decomposition 
for our space 
$c(K)$. For two simplicial complexes $K$ and $L$ we define the 
join product by $K*L\letbe \{(\sigma,\tau) : \sigma \in K, \tau \in L\}$. 
Then, 
by construction, we have $c(K*L)=c(K)\times c(L)$.  We also notice that, 
for the full simplex $\Delta(V)\letbe \{\sigma \subset V\}$ on a 
vertex set $V$,
we have $c(\Delta(V))\simeq BT^V$. This follows from the fact 
that $\catcat(\Delta(V))$ has $V$ as a terminal object and that 
therefore $\hoc_{\catcat(\Delta(V))} BT^{\Delta(V)} \simeq BT^K(V)=BT^V$.  
For every face $\sigma\in K$ we denote by 
$\st(\sigma)\letbe \st_K(\sigma)\letbe \{\tau\in K : \sigma\cup \tau \in K\}$ 
the star of $\sigma$.
This is again a simplicial complex and, 
since $\st(\sigma)=\Delta(\sigma)*\link(\sigma)$ 
we have 
$c(\st(\sigma))\simeq BT^\sigma\times c(\link(\sigma))$.Moreover, 
for $\sigma\subset \tau \in K$, we have $\st_K(\tau)\subset \st_K(\sigma)$ 
which induces a pointed map $c(\st_K(\tau)) \lra c(\st_K(\sigma))$.
This establishes  a functor
$$
\cst_K : \catkop \lra \Topx
$$
defined by $\cst_K(\sigma)\letbe c(\st_K(\sigma))$. Since the 
category $\catkop$ 
has an terminal object, namely $\emptyset$, 
we have an obvious homotopy equivalence
$\hoc^+_{\catkop} \cst_K \simeq \cst_K(\emptyset)=c(K)$. But 
restricting $\cst_K$ to the
full subcategory $\catkxop$ produces a map 
\nz
$\hoc^+_{\catkxop} \cst_K \lra c(K)$
which turns out to be an equivalence as well. 

\begin{thm}\label{seconddecomposition}$\ \ \ $
\nz
(i) $\hoc^+_{\catkxop}\ \cst_K \lra c(K)$ is a 
homotopy equivalence and a sharp $R$-homology decomposition.
\nz
(ii)
There exists a cofibration 
$$
|K| \lra \hoc_{\catkxop}\ \cst_K \lra c(K) .
$$
\nz
(iii) $$
\lim_{\catkxop}^{\ \ \ \ i} H^*(\cst_K;R)\cong
\begin{cases}
R(K)\oplus \widetilde H^0(K;R) & \text{for $i = 0$}\\
H^i(K;R) & \text{for $i>0$}.
\end{cases}
$$
\end{thm}

The rest of this section is devoted to a proof of this theorem.
We will compare the
two homotopy colimits,
$\hoc^+_{\catkxop} \cst_K$ and $\hoc^+_{\catk} BT^K$ 
and do this in several steps.
First we show that is sufficient to take the pointed homotopy colimit of
$BT^K$ over the category $\catkx$.

\begin{prop} \label{reductiontopointedlimit}
We have an equivalence
$$
\hoc^+_{\catkx} BT^K \larrow{\simeq} \hoc^+_{\catk} BT^K.
$$
Moreover, $\hoc^+_{\catkx} BT^K \lra c(K)$
is also  a sharp $R$-homology decomposition.
\end{prop}

\begin{proof} Since $\widetilde H^*(*)=0$, we have an isomorphism
$$
N(\catkxop, \widetilde H^*(BT^K))\cong N(\catkop, \widetilde H^*(BT^K))
$$
of normalised cochain complexes. This shows that the map between both pointed
homotopy colimits produces an 
isomorphism between 
higher limits, an isomorphism between the Bousfield-Kan
spectral sequences and therefore an isomorphism in cohomology. 
Moreover, both homotopy colimits are simply connected. Hence, 
by the Whitehead theorem, the map also induces an isomorphism
between the homotopy groups and is therefore a homotopy equivalence, 
which proves the second part.
\end{proof}

Now,
we construct a category $\catc$ which contains both, $\catkxop$ and  $\catkx$. 
This will allow to compare the pointed homotopy colimits
$\hoc^+_{\catkx} BT^K$ and
\nz
$\hoc^+_{\catkxop} \cst_K$.
To do this we will distinguish between the objects of
$\catkxop$, denoted by  $\tau\op$,
and the objects of $\catkx$, denoted by
$\tau$. The objects of $\catc$ are given by the union of the objects of
$\catkxop$ and $\catkx$. That is each face $\tau$ of $K$ generates
two objects in $\catc$, namely $\tau$ and $\tau\op$.
There exists at most one morphism between two objects.
And there are morphism
$\rho \ra \tau$, $\tau^{op} \ra \rho^{op}$ and $\tau \ra \rho^{op}$ 
if and only if $\rho \subset \tau$.
We have obvious inclusions $\phi: \catkxop \lra \catc$ and
$\psi:\catkx \lra \catc$.

Since for $\rho \subset \tau\in K$, the set $\tau$
is a face of $\st_K(\rho)$, there exists a well defined pointed map
$BT^\tau\lra c(\st_K(\rho))$. 
These maps are compatible with the pointed  inclusions
$BT^\alpha \subset BT^\beta$ as well as with
the pointed maps
$c(\st_K(\beta)) \ra c(\st_K(\alpha))$ for $\alpha \subset \beta$.
Therefore, we can define a functor $\bst: \catc \lra \Topx$
such that $\bst(\tau)\letbe BT^\tau$ and $\bst(\tau\op)\letbe \cst_K(\tau)$.
In particular, 
$\bst \psi=BT^K$ and $\bst \phi=\cst_K$.

Given a functor $\Phi : \catd' \lra  \catd$, for each object
$d\in D$, we can define the over category $\Phi\da d$. The objects
are given by morphisms $i':\Phi(d') \ra d$ in $D$, where $d'$
is an object of $\catd'$. And a morphism
$(i':\Phi(d') \ra d) \ra (i'':\Phi(d'') \ra d)$ is given by a morphism
$j:d' \ra d''$ of $\catd'$ such that $i''\Phi(j)=i'$.
The under category $d \da \rho$ is defined similarly.
As usual, we say that $\Phi$ is left cofinal if all over categories 
$\Phi\da d$ and right cofinal if all under categories
$d \da \rho$ are contractible; i.e. the classifying spaces are contractible.

The following series of statement shows how to compare the two 
homotopy colimits in question.

\begin{lem} \label{overundercategories} $\ \ \ $
\nz
(i)
For each object $\tau \in \catkx\subset \catc$, the under category
$\tau \da \phi$ is contractible.
\nz
(ii) For each object $\tau\op \in \catkxop \subset \catc$,
there exists an isomorphism of categories
$\catcat (\link_K(\tau)) \cong \psi \da \tau\op$ induced by
$\rho \mapsto (\rho\cup \tau \ra \tau\op)$.
\end{lem}

\begin{proof}
The first claim follows from the fact that $\tau \ra \tau\op$ is a
terminal object of
\nz
$\tau\da \phi$. The second claim follows
from an easy
straight forward calculation.
\end{proof}

\begin{prop} \label{firstcomparison} $\ \ \ \ $
\nz
(i) $\hoc_{\catkxop}^+ \cst_K \simeq  \hoc^+_{\catc} \bst$.
\nz
(ii) $\lim{}^i_{\catkx} \widetilde H^*(\cst_K) \cong
\lim{}^i_{\catc\op} \widetilde H^*(\bst)$
\end{prop}

\begin{proof}
Since for every object $\tau\op\in \catc$ the under category 
$\tau\op \da \phi$ is obviously contractible, Lemma \ref{overundercategories}
implies that  the inclusion
functor $\catkxop \ra \catc $ is right cofinal.
Since the restriction of $\bst |_{\catkxop} = \cst_K$, the equivalence
between the pointed homotopy colimits in
part (i) follows from \cite{boka}.

For the isomorphism in (ii) we need that
$(\catkxop)\op \ra \catc\op $ is left cofinal \cite{boka},
i.e. all under categories $c \da \phi\op \cong (\phi \da c)\op$ are
contractible. But this follows as above.
\end{proof}

\begin{prop} \label{secondcomparison} $\ \ \ $
\nz
(i) $\hoc_{\catkx}^+ BT^K \simeq  \hoc^+_{\catc} \bst$.
\nz
(ii) $\lim{}^i_{\catkxop} \widetilde H^*(BT^K)
\cong \lim{}^i_{\catc\op} \widetilde H^*(\bst)$.
In particular, for $i\geq 1$, 
\nz
$\lim{}^i_{\catc\op} \widetilde H^*(\bst)=0$.
\end{prop}

\begin{proof}
The left Kan extension $L\letbe L_{BT^K}$ of the functor 
$BT^K :\catkx \lra \Topx$
along the functor $\catkx \lra \catc$ is defined by
$L(c)\letbe  \hoc^+_{\psi \da c} BT^K$ . And
$\hoc_{\catc} L \simeq \hoc_{\catkx} BT^K$ \cite{boka}.
By Lemma \ref{overundercategories} and Theorem \ref{firstdecomposition},
there exists a natural transformation $L\lra \bst$, which induces a 
homotopy equivalence 
at each object. This proves the first part.

For the second part we apply the composition of functor spectral sequence
(e.g. see \cite{gazi} where the dual situation for colimits is discussed
in detail).
That is there exists a spectral sequence
$$
E_2^{i,j}\letbe \lim{}^i_{\catc\op} \lim{}^j_{c\da \psi^{op}} 
\widetilde H^*(BT^K) 
\Longrightarrow \lim{}^{i+j}_{\catkxop} \widetilde H^*(BT^K).
$$
By Theorem \ref{firstdecomposition} and Lemma \ref{overundercategories}
$\lim{}^j_{c\da \psi\op} \widetilde H^*(BT^K)=0$ for $j\geq 1$
and
$\lim{}_{c\da \psi\op} \widetilde H^*(BT^K)=\widetilde H^*(\bst(c))$.
This proves part (ii).
\end{proof}

\medskip

{\it Proof of Theorem \ref{seconddecomposition}.}
The first part follows from Proposition \ref{firstcomparison}
and Proposition \ref{secondcomparison}, and  the second part from the 
general relation between pointed and unpointed homotopy 
colimits as discussed above.

In the rest of the proof all limits are taken over $\catkxop$. 
Let $1_R$ denote the constant functor. Since a category and it's opposite 
category have the same geometric realization, 
we have $\lim{}^i 1_R\cong H^i(K)$.
The short exact sequence
$$
0 \lra \widetilde H^*(\cst_K) \lra H^*(\cst_K) \lra 1_R \lra 0 
$$
of functors establishes a long exact sequence of the higher limits.
By part (i), this long exact sequence splits into a short exact 
sequence 
$$
0 \lra \widetilde H^*(c(K))\cong \lim{}^0 \widetilde H^*(\cst_K) \lra 
\lim{}^0 H^*(\cst_K) \lra \lim{}^0 1_R\cong H^0(K) \lra 0
$$
and isomorphisms $\lim{}^i H^*(\cst_K)\cong \lim{}^i 1_R\cong H^i(K)$ for 
$i\geq 1$.
The short exact sequence can be rewritten as 
$$
0 \lra R(K) \lra \lim{}^0 R(\st_K) \lra \widetilde H^0(K)\lra 0,
$$
which proves part (iii).
\qed

\begin{rem} \label{atomicfunctor}
For later purpose we will calculate the higher limits for particular functors.
Let $M$ be an abelian group and $\psi_M:\catkop\lra \ab$ be the 
atomic functor
given by $\phi_M(\emptyset)\letbe M$ and $\phi_M(\sigma)\letbe 0$ for 
$\sigma\neq \emptyset$.Let $1_M : \catkop \lra \ab$ denote the 
constant functor 
which maps all objects to $M$ and all morphisms to the identity.
Then, we have a short exact sequence
$$
0\lra \phi_M \lra 1_M \lra \psi_M\letbe 1_M/\phi_M \lra 0.
$$
Since $\lim_{\catkop} 1_M\cong M\cong \lim_{\catkop} \psi_M$, the long
 exact sequence for the higher limits establishes isomorphisms
$\lim{}^i_{\catkop} \psi_M\cong \lim{}^{i+1} \phi_M$.
Since $N^*(\catkop;\psi_M)\cong N^*(\catkxop ; 1_M)$,
we get  a sequence of isomorphisms 
$$
H^i(K)\cong \lim{}^i_{\catkxop} 1_M \cong \lim{}^i_{\catkop} \psi_M \cong
\lim{}^{i+1}_{\catkop} \phi_M .
$$
By construction, this composition is natural with respect to maps 
between simplicial complexes. 
\end{rem}

\bigskip

\section{A vector bundle over $c(K)$} \label{vectorbundle}

In \cite{norac}, Theorem \ref{firstdecomposition} 
was used to construct a particular nontrivial map $c(K) \lra BU(n)$, whose  
construction we recall next.
 
Let $T^m \lra U(m)$ denote the maximal torus of the unitary 
group $U(m)$ given by 
diagonal matrices. 
The cohomology $H^*(BU(m);\Z)\cong \Z[c_1,...,c_m]$ of 
the classifying space $BU(m)$
is a polynomial algebra generated by the Chern classes $c_i$ and 
$H^*(BT^m;\Z)\cong \Z[V]$ is a polynomial algebra as well
which we identify with the polynomial algebra generated by the set
$V$ of vertices. The above map induces an 
isomorphism 
$H^*(BU(m);\Z)\cong \Z[V]^{\Sigma_m}\cong \Z[\sigma_1,...,\sigma_m]$,
where we identify $\Sigma_m$ with the Weyl group of $U(m)$ and 
where $\Sigma_m$ acts on $\Z[V]$ by permutations . Here, 
$\sigma_j$ denotes the $j$-th 
elementary symmetric polynomial. We can and will identify the Chern classes
$c_j$ with the elementary symmetric polynomials $\sigma_j$.

Since $\dim K=n-1$, a monomial $v_\tau$ vanishes in 
$\Z(K)$ if $\sharp\tau \geq n+1$. Hence, the composition 
$\Z[c_1,...,c_m] \lra \Z[V]\lra \Z(K)$ factors through 
$\z[c_1,...,c_n]$
and establishes a commutative diagram 
$$
\CD
\Z[c_1,...,c_m] @>>> \Z[V]\cong \Z[v_1,....,v_m]\\
@VVV @VVV \\
\Z[c_1,...,c_n] @>\psi>> \Z(K)
\endCD
$$
The left vertical arrow is induced by the canonical 
inclusion $BU(n) \lra BU(m)$, i.e. we put 
$c_j=0$ for $n+1\leq j\leq m$. The following statement 
concerning a topological 
realization of $\psi$ is proven in \cite{norac}.

\begin{thm} \cite{norac} \label{vectorbundel}
There exists a topological realization  $f_K :c(K) \lra BU(n)$ of $\psi$,
i.e. $H^*(f;\Z)=\psi$, which is unique up to homotopy. Moreover, 
the diagram
$$
\CD
c(K) @>f_K>> BU(n) \\
@VVV @VVV \\
BT^m @>>> BU(m)
\endCD
$$
commutes up to homotopy.
\end{thm}

The map $f_K$ is constructed as follows. Let $T^n\subset U(n)$ 
denote the maximal torus of $U(n)$ given by diagonal matrices. 
Since $\sharp \tau\leq n$, for each face $\tau \in K$, we can think 
of $\tau$ as a subset of the set $\{1,...,n\}$ which 
we also denote by $n$. Such 
an inclusion establishes a monomorphism 
$T^\tau \lra T^n \subset U(n)$ and, passing two classifying spaces, 
a map $f_\tau:BT^\tau \lra BU(n)$. Moreover, for a different inclusion 
$\tau\subset n$, the two monomorphisms $T^\tau \lra T^n$ 
differ only by a permutation. Hence, they are conjugate  in $U(n)$ 
and the induced maps  between the associated classifying spaces are homotopic
\cite{segal}.
This establishes a map from the 1-skeleton of the homotopy colimit
into $BU(n)$ unique up to homotopy. There exists an obstruction theory
for extending this map to a map $\hoc BT^K \lra BU(n)$ as well as 
for the uniqueness question  of such extensions \cite{wo}. The obstruction 
groups are higher limits of the form 
$\lim^j \pi_i(\map(BT^\tau,BU(n))_{f_\sigma})$.
If we pass to completion, i.e. we replace $BU(n)$ by it's p-adic completion
$BU(n)\p$,
the mapping space can be identified with 
$(BT^\tau\times BU(n\setminus \tau))\p$
\cite{no}.
For $j=i,i+1$ and target $BU(n)\p$, these higher limits 
do vanish \cite{norac}, 
which is sufficient to prove existence and homotopical uniqueness
of maps $f_K:c(K)\lra BU(n)\p$ realizing 
$\psi$ for all primes \cite{boka} \cite{wo}.
Rationally, the map $\psi$ can be realized, since the rationalisation 
$BU(n)_0\simeq \prod_{i=1}^n K(\Q,2i)$
of $BU(n)$ is a product of rational Eilenberg-MacLane spaces in even degrees.
An arithmetique square argument then establishes a map $f_K \lra BU(n)$ and
also shows that the homotopy class of this map is uniquely determined
(for details see \cite{norac}). 

As already mentioned in the introduction, we define $X_K$ as 
the homotopy fibre 
of $f_K : c(K) \lra BU(n)$. 
We are particularly interested in the top degrees of $H^*(X_K)$.

\begin{prop} \label{topdimensions} $\ \ \ \ $
\nz
(i)
$X_K$ has the homotopy type of a finite $CW$-complex of 
dimension $\leq n^2+n$.
\nz
(ii)
$$
H^i(X_K)\cong 
\left\{
\array{ll}
0 & \text{ for } i> n^2+n \\
 H^{n-1}(K) & \text{ for } i=n^2+n
\endarray
\right.
$$
\nz
(iii)
If $\widetilde H^i(K)=0$ for $i<n-1$, then 
$H^{n^2+n-1}(X_k)=0$ and there exists a short exact sequence
$$
0\lra H^{n^2+n-2}(X_K) \lra \prod_{i\in V} \widetilde H^{n-2}(\link_K(\{i\})) 
\lra H^{n-1}(K) \lra 0
$$
\end{prop}

\begin{proof}
In the proof cohomology is always
taken with coefficients in $R$. The composition 
$BT^\sigma \lra c(K) \larrow{f_K} BU(n)$ 
is natural with respect to maps in 
$\catk$. Interpreting this map as the classifying map of a 
$U(n)$-principal bundle
establishes a diagram of $U(n)$-principal bundles 
$Y(\sigma) \lra BT^\sigma$ 
over $\catk$. By construction, $Y(\sigma) \simeq U(n)/T^\sigma$. 
Since $U(n)$ acts freely on $Y(\sigma)$, the Borel construction 
$Y(\sigma)_{hU(n)}\letbe  Y(\sigma)\times_{U(n)} EU(n)$ 
projects to $Y(\sigma)/U(n)=BT^\sigma $
by a homotopy equivalence. These projections are natural with respect to
morphisms in $\catk$.
Since Borel constructions commute with taking homotopy colimits
\cite{drfa} we get a commutative diagram
of fibrations
$$
\CD
\hoc_{\catk} Y(-) @>>> \hoc_{\catk} BT^K @>>> BU(n) \\
@V\simeq VV @V\simeq VV @| \\
X_K @>>> c(K) @>>> BU(n)
\endCD
$$
We can calculate $H^*(X_K)$ with the help of the Bousfield-Kan 
spectral sequence.
Since the normalised cochain complex $N^i(\catkop, \phi)$ 
vanishes for $i>n$ for 
any functor $\phi$, we have 
$\lim{}^i_{\catk} \phi = 0$ for $i>n$.  
Moreover,
$H^i(U(n)/T^\sigma)=0$ for $i>n^2=\Dim U(n)$.
Since $\lim^i H^j(Y(-))$ is always a finitely generated abelian group,
this shows that $H^*(X_K)$ vanishes in degrees $> n^2+n$ and is a 
finitely generated abelian group in each degree. By construction,
$X_K$ is simply connected, and is therefore homotopy equivalent to a finite 
$CW$-complex of dimension $\leq n^2+n$.

The above argument also shows, that
the group $E_2^{n,n^2}\cong \lim{}^n_{\catk} H^{n^2}(Y(-))$ is the 
only possibly non vanishing 
term of total degree $n^2+n$ and survives in the spectral sequence.
In particular, $E_2^{n,n^2}\cong H^{n^2+n}(X_K)$. 
On the other hand, for any $\sigma\neq \emptyset$ we have
$\Dim U(n)/T^\sigma < n^2$. Hence,  
the functor
$H^{n^2}(Y(-))$ has it's only non vanishing value for $\sigma=\emptyset$
and $H^{n^2}(Y(\emptyset))\cong H^{n^2}(U(n))\cong R$.
Hence, by Remark \ref{atomicfunctor}, 
$E_2^{n,n^2}\cong \widetilde H^{n-1}(K)$. 
This proves 
the second part of the claim. 

In fact, Remark \ref{atomicfunctor} shows that 
$\widetilde H^{i-1}(K)\cong \lim{}^i_{\catkop} H^{n^2}(Y(-))\cong E_2^{i,n^2}$.
Hence, if $\widetilde H^j(K)$ vanishes for $j\neq n-1$, the only term of 
total degree $n^2+n-1$ is given by $\lim{}^{n} H^{n^2-1}(Y(-))$.
Let $\phi\letbe H^{n^2-1}(Y(-))$. 
Since $\Dim U(n)/T^\tau\leq n^2-\sharp\tau$, the functor $\phi$ vanishes for
all faces $\tau$ of order $\geq 2$. Moreover, for each vertex $i\in V$ 
the projection $U(n) \lra U(n)/T^{\{i\}}$ induces an isomorphism 
$H^{n^2-1}(U(n)/T^{\{i\}}) \larrow{\cong} H^{n^2-1}(U(n))$. This follows from
an analysis of the Serre spectral sequence of the fibration 
$T^{\{i\}}\lra U(n) \lra U(n)/T^{\{i\}}$.

For $r=0,1$ we define functors 
$\phi_r$ by $\phi_r(\sigma)\letbe \phi(\sigma)$ if 
$\sharp\sigma=r$ and $\phi_r(\sigma)\letbe 0$ otherwise.
We get short exact sequences of functors
$$
0 \lra \phi_0 \lra \phi \lra\phi_1 \lra 0
$$
and
$$
0 \lra \phi \lra 1_R \lra \psi\letbe 1_R/\phi \lra 0
$$
where $1_R$ denotes the constant functor.
The functor 
$\psi$ is non trivial only for faces of order $\geq 2$.
 
In  \cite{norab} higher limits of functors defined on $\catk$ 
are discussed 
in detail. Those results show, that $ \lim^j 1_R=0$ for $j\geq 1$, 
that $R=\lim{}^0 1_R \cong \lim{}^0 \phi$, that
$\lim{}^{j} \psi \cong \lim^{j+1} \phi$ for $j\geq 1$ 
and that $\lim^{j}\psi=0$ for 
$j\geq n-1$.
In particular,
$0=\lim{}^{n-1}\psi = \lim{}^n \phi$.
This proves that $H^{n^2+n-1}(X_K)=0$.

The first of the above sequences establishes an exact sequence
$$
0=\lim{}^{n-1} \phi_0 \lra \lim{}^{n-1} \phi \lra \lim{}^{n-1} \phi_1 
\lra \lim{}^n \phi_0 \lra \lim{}^n \phi =0
$$
By part (i) the fourth term can be identified with $\widetilde H^{n-1}(K)$, 
by \cite{norab}  the third term  with 
$\prod_{i\in V} \widetilde H^{n-2}(\link(\{i\}))$, and the second term with 
$H^{n^2+n-2}(X_K)$. This finishes the proof of the third part.
\end{proof}

\begin{cor} \label{restrictiontopdegree}
Let $L\subset K$ be a subcomplex of the same dimension.
Then, the composition
$$
H^{n-1}(K)\cong H^{n^2+n}(X_K) \lra H^{n^2+n}(X_L)\cong H^{n-1}(L)
$$
is the map induced in cohomology by the inclusion.
\end{cor}

\begin{proof}
This follows from the above proof and Remark \ref{atomicfunctor}
\end{proof}

We are also interested in the top degree of $H^*(X_{\st(\{i\})})$.

\begin{lem}
$X_{\st(\{i\})}$ has the homotopy type of a finite $CW$-complex 
of dimension $\leq n^2+n-2$ and 
$H^{n^2+n-2}(X_{\st(\{i\})})\cong H^{n-2}(\link(\{i\}))$.
\end{lem}

\begin{proof}
Since $c(\st(\{i\}))\simeq BT^{\{i\}}\times c(\link(\{i\}))$ 
we have a commutative
diagram of fibrations
$$
\CD
X_{\link(\{i\})} @>>> c(\link(\{i\}))\times BT^{\{i\}} 
@>>> BU(n-1)\times BT ^{\{i\}}\\
@VVV @V\simeq VV @VVV \\
X_{\st(\{i\})} @>>> c(\st(\{i\})) @>>> BU(n)
\endCD
$$
Since the homotopy fibre of the right vertical arrow is homotopy 
equivalent to the 
$n-1$-dimensional complex projective space $\CP(n-1)$,
this establishes a fibration 
$X_{\link(\{i\})} \lra X_{\st(\{i\})} \lra \CP(n-1)$.
This shows that $X_{\st(\{i\})}$ is simply connected and has 
the homotopy type of 
a finite $CW$-complex of dimension $\leq (n-1)^2+(n-1)+2(n-1)=n^2+n-2$ and 
that $H^{n^2+n-2}(X_{\st(\{i\})})\cong H^{(n-1)^2+(n-1)}(\link(\{i\}))$.
The last equation follows from an analysis of the 
Serre spectral sequence of the above
fibration.
\end{proof}

For later purpose we need the following lemma.

\begin{lem} \label{rkfinitegenerated}
$R(K)$ is a finitely generated $R[c_1,...,c_n]$-module.
\end{lem}

\begin{proof}
Since $R[V] \lra R(K)$ is a surjection and since $R[V]$ is a 
finitely generated $R[c_1,...,c_m]$-module, the same holds for
$R(K)$ as $R[c_1,...,c_m]$-module. By Theorem \ref{vectorbundel}
the map $R[c_1,...,c_m] \lra R(K)$ factors through 
$R[c_1,...,c_n]$, which implies the statement.
\end{proof}

We finish this section with the following observation:

\begin{rem} \label{n*structure}
The composition  $c(\st{(\tau)}) \lra c(K) \larrow{f_K} BU(n)$ makes 
$R(\st(\tau))$ into a $R[c_1,...,c_n]\cong H^*(BU(n);R)$-module and, 
with respect to this structure,
all differentials of the normalised chain complex 
$N^*(\catk ; H^*(c(\st_K);R)$ become
$H^*(BU(n);R)$-linear. 
Hence, $\lim{}^i H^*((\cst_K);R)=\lim{}^i R(\st_K)$ is an 
$H^*(BU(n);R)$-module.
The proof of Theorem \ref{seconddecomposition} shows that 
part (iii) can be refined. There exists a short exact sequence
$$
0\lra R(K) \cong H^*(c(K);R) \lra 
\lim{}^0 H^*(\cst_K);R) \cong \lim{}^0 R(\st_K) \lra \widetilde H^0(K;R) \lra 0
$$
of $H^*(BU(n);R)$-modules. Here, $H^*(BU(n);R)$ acts on 
$\lim{}^i R(\st_K)\cong H^i(K;R)$
as well as on $\widetilde H^0(K;R)$
via the augmentation $H^*(BU(n);R) \lra R$.
\end{rem}

\bigskip

\section{Cohen-Macaulay and Gorenstein conditions} \label{dcmg}

In this section we assume that $R=\bF$ is a field and cohomology is 
always taken with coefficients in $\bF$. In particular, 
$H^*(-)\letbe H^*(-;\bF)$.

Let $A^*$ be a non negatively graded commutative algebra over $\bF$.
We say that $A^*$ is connected if $A^0\cong \bF$ and 
$\bF$-finite if $A^j=0$ for $j$ large and
$A^j$ is a finitely generated $\bF$-module in each degree.
We call a finite sequence of
elements $a_1,..,a_r \in A$  a homogeneous system of
parameters, a hsop for short, if
they are homogeneous and algebraically independent and if
the quotient $A^*/(a_1,..., a_r)$ is  $\bF$-finite.
We  say that the sequence is a regular sequence,
if, for all $i$, $a_{i+1}$ is not a zero divisor in $A/(a_1,...,a_i)$.

We call $A^*$ Cohen-Macaulay, if there exists a sequence $a_1,...,a_n$ 
which is both,
a hsop and regular. If $A^*$ is Cohen-Macaulay, then it is known, that
every hsop is also a regular sequence \cite{brhe}.

A Noetherian local ring $S$ is called Gorenstein 
if $S$ considered as a module over itself 
has a finite injective resolution. If $A^*$ is a  commutative
connected non negatively graded  algebra, we can use the following 
equivalent definition \cite[Theorem I.12.4]{st}. That is,
$A^*$ is Gorenstein, if $A^*$ is Cohen-Macaulay
and if for any hsop $a_1,...a_n$ of $A^*$, we have 
$\soc(A^*/a_1,...,a_n)\cong \bF$. 
The socle $\soc(B^*)$ of  a non negatively graded algebra $B^*$ is defined as 
$\soc(B^*)\letbe \{b\in B^* : B^+b=0\}\cong \Hom_{B^*}(\bF,B^*)$ where
$B^+$ denotes the elements of positive degree. 

We call $A^*$ a Poincar\' e duality algebra, or PD-algebra for short,
if there exists a class $[A]\in \Hom_\bF(A^*,\bF)$ such that the composition
$A^*\otimes A^* \lra A^* \larrow{[A]} \bF$ is a non degenerate bilinear form. 
In particular,
every PD-algebra is connected.
As a straight forward argument shows, the above condition is equivalent to 
the fact that $\soc A^*\letbe \Hom_{A^*}(\bF, A^*)\cong \bF$. That is, 
a $\bF$-finite non negatively graded algebra is a $PD$-algebra 
if and only if it is Gorenstein. 
If $A^*\cong H^*(X;\bF)$ for a topological space, then we call $X$ a 
Poincar\'e duality space over $\bF$, a $\bF$-PD-space for short, if $A^*$ is a
PD-algebra. In this case $[X]\letbe [A]\in H_*(X)$ is the 
fundamental class of $X$.

We want to describe Cohen-Macaulay
and Gorenstein properties of the face ring $\bF(K)$ in terms of the map
$f_K: c(K) \lra BU(n)$ described 
in Theorem \ref{vectorbundel} and it's 
homotopy fibre $X_K$. The next result contains part (i) of 
Theorem \ref{charfibre}. 

\begin{thm} \label{cmcharfibre}
The following conditions are equivalent:
\nz
(i) $\bF(K)$ is Cohen-Macaulay.
\nz
(ii) The map $c(K) \lra BU(n)$ makes $\bF(K)$ into a finitely generated
free $H^*(BU(n))$-module.
\nz
(iii) The cohomology ring  $H^*(X_K)$ is concentrated
in even degrees.
\nz
(iv) $\Tor^j_{H^*(BU(n))}(\bF(K),\bF)=0$ for $j\leq -1$.
\end{thm}

For the proof  we need the following lemma. Again we 
denote by $\sigma_j$ the 
$j$-th elementar symmetric polynomial in the generators of $\bF(K)$ and 
identify $H^*(BU(n))$ with $\bF[\sigma_1,...,\sigma_n]$, 
that is with image of the
map $H^*(BU(n)) \lra H^*(c(K))$.

\begin{lem} \label{regularsequence}
The sequence $\sigma_1,..,\sigma_n\in \bF(K)$ is a hsop for $\bF(K)$.
\end{lem}

\begin{proof}
By Lemma \ref{rkfinitegenerated}
$\bF(K)/(\sigma_1,...,\sigma_n)$ is $\bF$-finite.
We have only to show that the elements $\sigma_1,...,\sigma_n$ are 
algebraic independent 
in $\bF(K)$.

Let $\mu$ be a maximal face of $K$, that is $\sharp \mu=n$.
The composition $BT^\mu \lra c(K) \lra BU(n)$ is induced
by a maximal torus inclusion $T^\mu \lra U(n)$ 
(see Section \ref{vectorbundle}). 
The images of $\sigma_1, ..., \sigma_n$ in $H^*(BT^\mu)\cong \bF[\mu]$ are 
given by the 
elementary symmetric polynomials and therefore 
are algebraic independent as well as $\sigma_1,...\sigma_1$ in $\bF(K)$.
\end{proof}

\medskip

{\it Proof of Theorem \ref{cmcharfibre}:}
If $\bF(K)$ is Cohen-Macaulay, then every hsop is given by a regular
sequence. And if we have a hsop of $\bF(K)$ given by a
regular sequence, then $\bF(K)$ is Cohen-Macaulay  \cite{brhe}.
In the light of Lemma \ref{regularsequence} this shows that the
first two conditions are equivalent.

If $H^*(X_K;\z)$ is concentrated in even degrees (part (iii)),
then, by degree
reasons, the Serre spectral sequence for the fibration
$X_K \lra c(K) \lra BU(n)$ collapses at the $E_2$ page
and $H^*(c(K)) \cong H^*(X_K)\otimes_\bF H^*(BU(n))$
as $H^*(BU(n))$-module. This shows that
$\bF(K)$ is a finitely generated free $H^*(BU(n))$-module
and therefore
Cohen-Macaulay, which is part (i).

If $\bF(K)$ is Cohen-Macaulay, then, by (ii),
it is a finitely generated free module over $H^*(BU(n))$.
In particular, $\Tor^j_{H^*(BU(n))}(\bF(K),\bF)=0$ for $j\leq -1$. This is
part (iv).

If  condition (iv) is satisfied, the Eilenberg-Moore spectral sequence
for calculating 
\nz
$H^*(X_K;\bF)$
collapses at the $E_2$-page and shows that
$H^*(X_K;\bF) \cong \bF(K) \otimes_{H^*(BU(n);\bF)} \bF$ and that 
$H^*(X_K;\bF)$ is
concentrated in even degrees, which is condition (iii). This proves the
equivalence of  (i), (iii) and (iv).
\qed

\medskip

{\it Proof of Theorem \ref{charfibre} (ii):}
By our definition of Gorenstein,
$\bF(K)$ is Gorenstein if and only if $\bF(K)$ is Cohen-Macaulay 
and
$\bF(K)\otimes_{\bF[\sigma_1,...,\sigma_n]} \bF 
\cong \bF(K)/(\sigma_1,...,\sigma_n)$  is
a PD-algebra. Hence the equivalence of the two conditions follows from 
the first part of the 
Theorem. 
\qed

\bigskip

\section{Homotopy fixed point sets} \label{hofipo}

For an action of a group $G$ on a space $X$ we can think of 
the fixed point set as the mapping space 
$X^G=\map_G(*,X)$ of $G$-equivariant maps from a point into $X$.
The notion is not flexible enough for doing homotopy theory, since a
homotopy equivalence between $G$-spaces, which happens to be 
$G$-equivariant in addition, does not induce an equivalence between 
the fixed-point sets in general. We therefore are 
interested in homotopy fixed point sets which do have this property. 
They are defined as 
the equivariant mapping space $X^{hG}\letbe \map_G(EG,X)$ 
where $EG$ is a contractible 
$G$-CW-complex with a free $G$-action. The projection $EG \lra *$
induces a map
$X^G \lra X^{hG}$.

Applying the Borel construction establishes a fibration
$$
X \lra X_{hG}\letbe X\times_G EG \larrow{\pi} BG. 
$$
A straight forward argument 
shows that  we can equivalently
 define
the homotopy fixed point set as the space 
$\Gamma(X_{hG} \ra BG)$ of sections
of this fibration.
The latter definition also allows to define homotopy fixed point sets 
in more general situations. A proxy $G$-action on $X$ is a fibration
$X \lra E \larrow{\pi} BG$, where we think of $E$ as the Borel construction
of this action. We define $X^{hG}\letbe \Gamma(E \ra BG)$. This establishes 
a fibration
$$
X^{hG} \lra \map(BG,E)_{\{id\}} \lra \map(BG,BG)_{id}
$$
Here, the middle term consist of all lifts of the identity $id$ of $BG$ 
up to homotopy, 
i.e. of all
maps $g:BG\lra E$ such that $\pi f\simeq id$. If $G$ is a 
finite abelian group the base space is 
homotopy equivalent to $BG$ \cite{hu}, and composition of maps 
yields an equivalence $BG\times X^{hG} \larrow{\simeq} \map(BG,E)_{\{id\}}$.
This equivalence fits into a commutative diagram of fibrations
$$
\CD
X^{hG} @>>> BG\times X^{hG} @>>> BG \\
@VVV @VVV @| \\
X @>>> E @>>> BG
\endCD
$$
where the vertical arrows are induced by evaluation at a basepoint.

Typical examples of proxy actions arise from pull back constructions.
Let 
\nz
$X \lra E \larrow{\pi} B$ be a fibration and $f: BG \lra B$ be a map.
The pull back construction establishes a commutative diagram
$$
\CD
E' @>\pi'>> BG \\
@VVV @VfVV \\
E @>\pi>> B
\endCD
$$
and applying the mapping space functor yields a pull back diagram
$$
\CD
\map(BG,E')_{\{id\}} @>\pi'_*>> \map(BG,BG)_{id} \\
@VVV @VfVV \\
\map(BG,E)_{\{f\}} @>\pi_*>> \map(BG,B)_{f}.
\endCD
$$
The fibre of both horizontal arrows is given by the 
homotopy fixed point set $X^{hG}$.
Here, the  left mapping space in the bottom row 
consists of all maps $BG\lra E$ which are
homotopic to a lift of $f$.
If $G$ is finite and abelian
the composition $BG\simeq \map(BG,BG)_{id} \lra \map(BG,B)_f \larrow{ev} B$ 
equals the map $f : BG \lra B$. 

The main goal of this section is to show that, in favourable cases,
$H^*(X^{hG};\fp)$ 
is concentrated in even degrees or a PD-algebra,  if 
$H^*(X;\fp)$ satisfies these properties. 
For our method of proof we 
have to make some restrictions. According to our situation, 
we have the spaces $X_K$ in mind, 
we will always assume that:
\begin{enumerate}
\item $X$ is $\fp$-finite and p-complete.
\item $H^*(X;\fp)$ is concentrated in even degrees.
\end{enumerate}
In particular, $H^*(X;\fp)$ is a graded algebra, commutative in the
non graded sense. Spaces satisfying both conditions will be called special.

We also restrict ourselves to particular proxy actions.
We say that a proxy action $X \lra E \lra BG$ is orientable if
$\pi_1(BG)$ acts trivially on $H^*(X;\Z\p)$.
 
\begin{rem} \label{actionassumption}
If $G$ is an elementary abelian p-group, and the $G$-action extends to a 
torus action of $T^r$, i.e. the proxy action is induced by a pull back of a
fibration $X \lra E' \lra BT^r$, 
the group $G$ always acts trivially on the cohomology of $X$ and the 
proxy action is orientable.

Moreover, if the proxy action $X\lra E \lra BG$ is orientable, the
Serre spectral 
sequence  for $H^*(-;\z\p)$ collapses
at the $E_2$-page by degree reasons. And the same holds for $H^*(-;\fp)$. 
In particular, $H^*(E;\fp)\cong H^*(BG;\fp)\otimes H^*(X;\fp)$
as $H^*(BG;\fp)$-module \cite[Corollary 2.5]{dela} . 
\end{rem}

The first part of the next theorem is due to Dehon and Lannes
\cite[Corollary 2.10]{dela}.

\begin{thm} \label{hofixinherited}
Let $G$ be an elementary abelian $p$-group and $X \lra E \lra BG$
an orientable  proxy $G$-action on a space $X$.
\nz
(i) If $X$ is special, then so is $X^{hG}$.
\nz
(ii) If $X$ is special and a $\fp$-PD-space, then so is each component of 
$X^{hG}$.
\end{thm}
 
The proof of the second part needs some preparation. For the rest 
of this section we 
make the following assumptions.
$H^*(-)$ denotes $H^*(-;\fp)$, $G$ is an elementary abelian $p$-group
and $H^G\letbe H^*(BG)$. Moreover, 
$X$ is special and  the proxy $G$-action $X \lra E \lra BG$ is 
orientable. In particular, $H^*(E)\cong H^G\otimes H^*(X)$ as $H^G$-module
(see Remark \ref{actionassumption}). 

For such actions $H^*(E)$ is not concentrated in even degrees
and hence not commutative in the non graded sense.
To avoid technical difficulties we will use the following construction.
We denote  by $J\subset H^G$ the
ideal generated by all classes of degree 1. And for 
an graded $H^G$-module $M$ 
we denote by $\bar M$ the quotient $M/JM$. 
If $G\cong (\Z/p)^r$,
the composition $H^*(BT^r) \lra H^G \lra \bar H^G$ is an 
isomorphism, where the first map is induced by the canonical inclusion 
$G\cong (\Z/p)^r \subset T^r$. In particular, 
$\bar H^G$ is a polynomial algebra 
generated by elements of degree 2.  In fact, we can think of it as the polynomial 
part and as a subalgebra
of
$H^G$. Hence, every $H^G$-module is naturally an $\bar H^G$-module.

\begin{lem} \label{hatgorenstein}
$H^*(X)$ is a PD-algebra if and only if $\bar H^*(E)$ is Gorenstein.
\end{lem}

\begin{proof}
By assumption, $H^*(E)\cong H^*(X)\otimes H^G$ as $H^G$-module.
And hence, $\bar H^*(E)\cong H^*(X)\otimes \bar H^G$ as 
$\bar H^G$-module. Since $\bar H^G$ is a polynomial algebra,
this implies that $\bar H^*(E)$ is Cohen-Macaulay. And hence,
$\bar H^*(E)$ is Gorenstein if and only if 
$H^*(X)\cong \bar H^*(E)\otimes_{\bar H^G} \fp$ is a PD-algebra.
\end{proof}

 Let 
$S\subset H^G$ denote the multiplicative subset generated by all 
Bockstein images in  degree 2 of nontrivial elements of degree 1. 
That is the subset of all images of non trivial elements of $H^*(BT^r)$ 
of strictly positive degree. For any 
$H^G$-module $M$ we denote by $S^{-1}M$ the localised module over 
$S^{-1}H^G$.
Let $X\lra E \lra BG$ be a $G$-proxy action on $Y$.
Following \cite{dwwi}
(Corollary 1.2 and the following remark), there exists a map
$S^{-1}H^*(E) \lra S^{-1}(H^G\otimes H^*(Y^{hG}))$ 
between the localised modules. Under favourable 
circumstances this is an isomorphism, which, as 
shown in \cite[Section 2]{dela}, always hold if
$X$ is special and if the $G$-action is orientable.
We collect this into the following theorem.

\begin{thm} (\cite{dwwi} \cite{dela}) \label{smiththeory}
Let $X \lra E \lra BG$ be an orientable $G$-proxy action on a special space $X$.
Then there exist isomorphisms 
$$
S^{-1}H^*(E) \cong S^{-1}(H^G\otimes H^*(X^{hG}))\cong 
(S^{-1}H^G)\otimes H^*(X^{hG}).
$$
\end{thm}

Since the multiplicative subset $S\subset H^G$ consist of elements of even degree, 
it also is  a multiplicative subset of  $\bar H^G$. 
We can think of the localised module 
$S^{-1}\bar H^*(E)\cong 
\bar H^*(E)\otimes_{\bar H^G} S^{-1}\bar H^G$ in a different way.
Since $\bar H^G\lra \bar H^*(E)$ is a monomorphisms, 
the set $S$ gives rise to a
multiplicative 
subset $S'\subset \bar H^*(E)$. Then multiplication induces an isomorphism
$\bar H^*(E)\otimes_{\bar H^G} S^{-1}\bar H^G \cong 
(S')^{-1} \bar H^*(E)$,
where the latter localisation is obtained by localising the algebra
$\bar H^*(E)$ with respect to the subset $S'$.

\medskip

{\it Proof of Theorem \ref{hofixinherited} (ii):}
 Let us assume that $X$ is an $\fp$-PD-space. Then, $\bar H^*(E)$ is 
Gorenstein (Lemma \ref{hatgorenstein}) as well as $S^{-1}\bar H^*(E)$
\cite[Proposition 3.1.19]{brhe}.
Since 
$$
S^{-1}\bar H^*(E)\cong H^*(X^{hG})\otimes S^{-1}\bar H^G \cong 
\bigoplus_g H^*(X^{hG}_g) \otimes S^{-1}\bar H^G
$$
as algebras, each of the summands is Gorenstein.
Here, the direct sum is taken over  components of  $X^{hG}$. Let 
$X^{hG}_g$ denotes the component associated to a section $g: BG \lra E$ of 
$\pi : E \lra BG$. Now we give the algebra 
$H^*(X^{hG}_g) \otimes S^{-1}\bar H^G$ 
a different grading induced by the grading of the first factor. 
That is elements of 
$S^{-1}H^G$ get degree 0.
Since $\bar H^G\cong \fp[t_1,...,t_n]$, the sequence 
$\{t_1-1,...t_n-1\}\subset H^*(X^{hG}_g) \otimes S^{-1}\bar H^G$ is  regular,
homogeneous and consists of elements of degree 0. Hence 
$H^*(X^{hG}_g) \cong 
H^*(X^{hG}_g) \otimes S^{-1}\bar H^G /(t_1-1,...,t_n-1)$
is a connected $\fp$-finite Gorenstein algebra and therefore a 
PD-algebra (see Section \ref{dcmg}). 
\qed

\bigskip

\section{Cohen-Macaulay and Gorenstein properties for links of faces}

In this section we want to prove that  Cohen-Macaulay or Gorenstein 
properties are inherited to the Stanley Reisner algebras of links of faces.
We call a simplicial complex pure, if all maximal faces have 
the same dimension.
We first consider the case of algebras over $\fp$.

\begin{thm} \label{proplinkfp}
If $\fp(K)$ is Cohen-Macaulay respectively Gorenstein, 
then $K$ is pure and, for each face $\tau\in K$
the algebra $\fp(\link_K(\tau))$ is also Cohen-Macaulay 
respectively Gorenstein.
\end{thm}

For the proof we will use the space $X_K$ given by the fibration 
$X_K \lra c(K) \lra BU(n)$ described in Section \ref{vectorbundle}.
Since $BU(n)$ is simply connected, the p-adic completion maintains 
the fibration
\cite{boka}.
Since we are working  with $\fp$ as coefficients and since for 
simply connected 
spaces
completion induces an isomorphism in mod-p 
cohomology, we can and will assume that all 
spaces are completed. To simplify notation, we will always drop the notation 
for completion.
The proof also relies on the interpretation of $\fp(\link(\tau))$ as a
certain mapping space
given in \cite{norab}, which we recall next.

For each face $\tau\in K$ we denote by $G^\tau \subset T^\tau$
the maximal elementary abelian subgroup of the torus $T^\tau$.
Let $g_\tau : BG^\tau \lra c(K)$ denote the composition 
$BG^\tau \lra BT^\tau \lra \hoc_{\catk} BT^K \simeq c(K)$.
The composition 
$BG^\tau \larrow{g_\tau} c(K) \larrow{f_K} BU(n)$ 
establishes a proxy $G^\tau$-action
$X_K \lra E \lra BG^\tau$. 
Applying the mapping space functor we get a fibration
$$
X^{hG^\tau} \lra \map(BG^\tau,c(K))_{\{f_Kg_\tau\}} \larrow{(f_K)_*} 
\map(BG^\tau,BU(n))_{f_Kg_\tau}
$$
(see Section \ref{hofipo}). 
The composition $f_Kg_\tau$ is induced from 
a coordinate-wise inclusion $G^\tau \subset T^\tau \subset T^n\subset U(n)$ 
into the set of diagonal matrices (see Section \ref{vectorbundle}).
The centraliser of this image equals 
$T^\tau\times U(n\setminus \tau)$ where we again think 
of $n$ as the set $\{1,...,n\}$ and $\tau$ becomes a subset of the set
$n$ via the coordinate-wise inclusion. Then, by construction, $f_Kg_\tau$ 
factors through 
a map 
$\id\times\const: BT^\tau \lra BT^\tau \times BU(n\setminus\tau)$ 
where the fist coordinate is the identity
and the second the constant map. There also exists a map
$BT^\tau\times c(\link(\tau))=c(\st(\tau)) \lra c(K)$ 
as constructed in Section \ref{hodec}
and the map $g_\tau$ factors through 
$i\times \const : BG^\tau \lra BT^\tau \times c(\link(\tau))$.
Moreover all these maps fit into a diagram
$$
\CD
X_{\link(\tau)} @>>> BT^\tau \times c(\link(\tau)) 
@>id\times f_{\link(\tau)}>>
BT^\tau \times BU(n\setminus \tau) \\
@VVV @VVV @VVV \\
X_K @>>> c(K) @>f_K>> BU(n).
\endCD
$$
Applying the mapping space functor, there exit  the following equivalences 
(of p-completed) spaces;
$BT^\tau \larrow{\simeq} \map(BG^\tau,BT^\tau)_{i}$, 
$BT^\tau\times BU(n\setminus\tau) \larrow{\simeq} 
\map(BG^\tau, BU(n)_{f_Kg_\tau}$
\cite{dwza} 
and
$BT^\tau\times c(\link(\tau)) \larrow{\simeq} \map(BG^\tau,c(K))_{g_\tau}$
\cite{norab}.
Putting all this information together 
we have a commutative diagram
$$
\CD
X_{\link(\tau))} @>>> BT^\tau\times c(\link(\tau)) @>f_{\link(\tau))}>> 
BU(n\setminus\tau) \\
@VVV  @V\simeq VV  @V\simeq VV \\
(X_K)^{hG^\tau}_{g_\tau} @>>> \map(BG^\tau,c(K))_{g_\tau} @>>> 
\map(BG^\tau,BU(n))_{f_Kg_\tau}
\endCD
$$
of horizontal fibrations.
Since $BT^\tau\times BU(n\setminus\tau)$ is simply connected,
the map 
$X_K^{hG^\tau} \lra \map(BG^\tau , c(K))_{\{f_Kg_\tau\}}$ 
induces a bijection between the
components of both spaces. Hence, the bottom
left space is connected.
This proves the following proposition.

\begin{prop}
$X_{\link(\tau)}\simeq (X_K)^{hG^\tau}_{g_\tau}$.
\end{prop}

Now we are in the position to prove Theorem \ref{proplinkfp}.

\medskip

{\it Proof of Theorem \ref{proplinkfp}}
If $\fp(K)$ is Cohen-Macaulay, the fibre $X_K$ is special 
(Theorem \ref{charfibre}(i)).
Since the map $f_Kg_\tau : BG^\tau \lra BU(n)$ factors through $BT^\tau$, 
the proxy $G^\tau$ action extends to a torus action and is orientable
(see Remark \ref{actionassumption}). 
We can apply
Theorem \ref{hofixinherited}. 
That is that $(X_K)^{hG^\tau}_{g_\tau}\simeq X_{\link(\tau)}$ 
is special and 
that $\fp(\link(\tau))$ is Cohen-Macaulay (Theorem \ref{charfibre}(i)).
If $\fp(K)$ is Gorenstein, we use Theorem \ref{charfibre}(ii) instead
of the first part.

Finally we have to show that $K$ is pure. 
The above argument shows that, if $\fp(K)$ is
Cohen-Macaulay, then $\fp(\link(\tau))$ is a free 
$H^*(BU(n\setminus \tau))$-module.
If $\mu \in K$ is a maximal simplex, then $\link(\mu)$ is 
the empty complex and
$\fp(\link(\mu))\cong \fp$. This implies that $\mu$ has 
order $n$ and that $K$ is pure.
\qed

\smallskip

We finally give up the restriction on the coefficients.

\begin{cor} \label{proplink}
If $\bF(K)$ is Cohen-Macaulay respectively Gorenstein, then $K$ is pure and, 
for every face 
$\tau \in K$, the Stanley-Reisner algebra 
$\bF(\link_K(\tau))$ is also Cohen-Macaulay 
respectively Gorenstein. 
\end{cor}

\begin{proof}
We will use Theorem \ref{proplinkfp} and Theorem \ref{charfibre}.
If $\bF$ is a field of characteristique $p>0$, then 
$H^*(X_K)\cong H^*(X_K;\fp)\otimes_{\fp} \bF$. This shows that
$\bF(K)$ is Cohen-Macaulay if and only if $\fp(K)$ is so. The same holds for 
links. This covers the Cohen-Macaulay case as well 
as the Gorenstein case and also shows that $K$ is pure.

Let $\z_{(p)}$ denote the localisation of $\z$ at the prime $p$.
Then, for $\bF=\Q$ we can argue as follows. Since $X_K$ is of the 
homotopy type of a finite $CW$-complex (Proposition \ref{topdimensions}),
$H^*(X_K;\Z)$ has only finitely many torsion primes. Hence, 
$\Q(K)$ is Cohen-Macaulay if and only if 
$H^*(X_K;\q)$ is concentrated in even degrees if and only if, 
for almost all  primes, 
$H^*(X_K;\z_{(p)})$ is 
concentrated in even degrees and torsion free
if and only if, 
for almost all  primes,
$H^*(X_K;\fp)$ is concentrated in even degrees
if and only if, 
for almost all  primes,
$\fp(K)$ is Cohen-Macaulay. 
Now, Theorem \ref{proplinkfp} and and the above chain of 
equivalent statements
applied in the case of $\link(\tau)$ proves the claim for Cohen-Macaulay 
algebras over $\bF=\q$.
In the Gorenstein case, we only have to notice that $H^*(X_K;\Q)$ is 
a $PD$-algebra if and only $H^*(X_K;\fp)$ is a $PD$-algebra for 
almost all primes.

For a general field of characteristique 0 we deduce the 
claim from the case $\bF=\Q$ 
in the same manner as for fields of characteristique $p>0$ from $\bF=\fp$.
\end{proof}

\bigskip

\section{Proof of Theorem \ref{cmresult}}
 
Theorem \ref {cmresult} follows  
easily from the following statement by induction.

\begin{thm} \label{cmequivalentstatements}
$\bF(K)$ is Cohen-Macaulay if and only if  $\bF(\link_K(\tau))$ is 
Cohen-Macaulay for 
all faces $\tau\neq \emptyset$ of $K$ and 
$\widetilde H^r(K;\bF)=0$ for $0\leq r<n-1$.
\end{thm}

\medskip

{\it Proof of Theorem \ref{cmresult}}:
If $\Dim K=0$ then $\bF(K)$ is always Cohen-Macaulay. In fact,
in this case $c(K)$ is the m-fold wedge product of $BS^1$'s and $X_k$ the 
$(m-1)$-fold wedge product of $S^2$'s, whose cohomology is concentrated in 
even degrees. On the other hand the set of conditions on the cohomology 
of the links of $K$ is an empty set. This proves the statement in this case.

The general case follows by induction over 
the dimension of $K$ and the above  
theorem.
\qed

\medskip

In the following cohomology is always  with $\bF$-coefficients. We define  
$H^*(-)\letbe H^*(-;\bF)$
for the rest of this
section. Also, for simplification, we set 
$P\letbe H^*(BU(n))\cong\bF[\sigma_1,...,\sigma_n]$.  

 The proof of Theorem \ref{cmequivalentstatements} 
is based on the homological analysis of  particular double complexes. 
To fix notation, we will recall the general concept next. For details see
\cite{mccleary}
 
Let $A$ be a $\bF$-algebra. A differential graded $A$-module  $(C^*,d_C)$ 
is a cochain complex 
of $A$-modules such that $d_C$ is $A$-linear. 
A double complex  or differential 
bigraded module $(M^{*,*}, d_h,d_v)$ over $A$
is a bigraded $A$- module $M^{*,*}$
with two $A$-linear maps
$d_h : M^{*,*} \lra M^{*,*}$ and $ d_v:M^{*,*} \lra M^{*,*}$ of bidegree
$(1,0)$ and $(0,1)$ such that $d_hd_h=0=d_vd_v$ and $d_hd_v+d_vd_h=0$.
We think of $d_h$ as the horizontal and of $d_v$ as
the vertical differential.
To each double complex $M^{*,*}$ we associate a total complex $Tot^*(M)$
which is  a differential graded module over $A$ defined by
$Tot^n(M)\letbe \oplus_{i+j=n} M^{i,j}$ with differential $D\letbe d_h + d_v$.
Examples are given by the tensor products of two differential graded
modules $(B^*,d_B)$ and $(C^*,d_C)$. If we set
$M^{i,j}\letbe B^i\otimes_A C^j$, $d_h\letbe  d_B\otimes 1$ and
$d_v\letbe (-1)^j 1\otimes d_C$,
we get a double complex such that $Tot^*(M)=B^*\otimes_A C^*$.

For a double complex $(M^{*,*}, d_h,d_v)$, we can take horizontal or
vertical cohomology groups denoted by 
$H^*_h(M^{*,*})$ and $H^*_v(M^{*,*})$.
The boundary maps $d_h$ and $d_v$ induce again boundary maps on
these cohomology groups. We can consider cohomology groups of the form
$H^*_h(H^*_v(M^{*,*}))$ and $H^*_v (H^*_h(M^{*,*}))$.

If $(M^{*,*}, d_h,d_v)$ is bounded below, that is $M^{i,j}=0$ if $i$ or
$j$ is small enough,
there exist two
spectral sequences converging towards $H^*(Tot(M),D)$.
In one case, we have
$E_2^{i,j}=H^i_h(H^j_v(M))$ and in the other case
$E_2^{i,j}\letbe H_v^j(H^i_h(M))$. 
In the first case the differential have degree $(r,1-r,)$ and in
the second case degree $(1-r,r)$. 

Let $N^*\letbe N^*(\catkxop ,H^*(\cst_K))$ denote 
the normalised cochain complex
for the functor $H^*(\cst_K) : \catkxop \lra \Ab$ 
considered in Section \ref{hodec}.
Actually, the complex $N^*$ is a bigraded object. 
It inherits an internal degree
from the grading of $H^*(\cst_K)$, which we will not consider in 
most cases. 
 
We collect the main properties of $N^*$ in the next proposition.

\begin{prop} \label{propertyn} $\ \ \ $
\nz
(i) $N^i=0$ for $i<0$ or $i\geq n$.
\nz
(ii) $(N^*,d_N)$ is a differential graded $P$-module as well as 
$H^*(N^*,d_N)\cong \lim\!{}_{\catkxop}^*\ , \bF(\st_K)$
\nz
(iii) There exist an isomorphism
$H^i(N^*,d_N) \cong \widetilde H^{i}(K)$ for $i\geq 1$ of $P$-modules and 
a short exact sequence  
$$
0\lra \bF(K) \lra H^0(N^*,d_L) \cong \lim\!{}_{\catkxop}{}^0 \,\bF(\st_K) \lra 
\widetilde H^0(K)\lra 0
$$
of $P$-modules.
\end{prop}

\begin{proof}
The first part follows from the fact that $N^s({\catkxop},H^*(\cst_k))=0$ 
for $s\geq n$,
the second and the third part from Remark \ref{n*structure}.
\end{proof}

We can also look at the projective resolution of the trivial $P$-module $\bF$ 
given by the 
Koszul complex $Q^*\letbe \Lambda^*\otimes_\bF P$. According to our 
convention we make this into a 
cochain complex and give the generators of the exterior algebra 
$\Lambda^*\letbe \Lambda^*(x_1,...,x_n)$ the degree $-1$. 
 As usual the differential 
$d_Q$ is defined by $d_Q(x_i)\letbe \sigma_i$ and $d_Q(y)=0$ for $y\in P$ 
and has 
degree 1.
Again, $Q^*$ is  a differential graded $P$-module, bounded below 
and above by $Q^j=0$ for $j>0$ or $j<-n$. 
The differential bigraded $P$-module
$N^*\otimes_{P} Q^*$ is then bounded.
In particular, both above mentioned 
spectral sequences converge towards
$H^*(\Tot(N^*\otimes_PQ^*))$. 

\medskip

{\it Proof of Theorem \ref{cmequivalentstatements}:}
If $K$ is the empty complex, there is nothing to show. If $K$ is a 
0-dimensional complex, then $\bF(K)$ is always 
Cohen-Macaulay as discussed in
the proof of Theorem \ref{cmresult} and Theorem \ref{gorensteinresult}. 

Now we assume that $\Dim K\geq 1$. i.e. $n\geq 2$.
We start with the assumption that $\bF(K)$ is Cohen-Macaulay.
By Theorem
\ref{proplink}, we know that $K$ is
pure and that for all
$\tau \in K$, the algebra $\bF(\link(\tau))$ is also
Cohen-Macaulay as well as
$\bF(\st(\tau))$.
In particular, since $\Dim \st(\tau)=\Dim K$, the algebra
$R(\st_K(\tau))$ is a finitely generated free module over $P$.
We have to show that $\widetilde H^i(K)=0$ for $i<n-1$.

We look at the above constructed double complex $N^*\otimes_{P} Q^*$.
All modules $N^i$ and $Q^j$ are free $P$-modules.
Hence, the functors $N^*\otimes_{P} -$ and $-\otimes_{P} Q^*$ are exact.
We get
$$
H^i_h (H^j_v (N^*\otimes_{P} Q^*)) \cong 
H_h^i(N^*\otimes_{P} H^j_v(Q^*))\cong
\left\{
\begin{array}{ll}
0 & \text{ for } j \neq 0 \\
H^i(N^*\otimes_{P} \bF) & \text{ for } j=0
\end{array}
\right.
$$
In particular, the $E_2$-term is concentrated in one  horizontal line
given by $j=0$,
the spectral sequence collapses and
$H^i(\Tot(N^*\otimes_{P} Q^*))\cong 
H^i(N^*\otimes_{P} \bF)=0$ for $i< 0$
or $i\geq n$.

Considering the second spectral sequence
we get
$$
\begin{array}{l}
H_v^j(H^i_h(N^*\otimes_{P} Q^*))  \cong 
H_v^j(H^i_h(N^*)\otimes_{P} Q^*) \\
 \hspace*{1cm} \cong
\left\{
\begin{array}{ll}
\Tor^j_{P}(\widetilde H^{i}(K),\bF) \cong 
\Tor^j_{P}(\bF,\bF)\otimes_\bF \widetilde H^{i}(K) & \text{ for } i > 0 \\
H^*(X_K)\oplus \Tor^0_P(\bF,\bF)\otimes \widetilde H^0(K) 
& \text{ for } i=0 \text{ and } j=0 \\
\Tor^0_P(\bF,\bF)\otimes \widetilde H^0(K) & 
   \text{ for } i=0 \text{ and } j\neq 0 \\
0 & \text{ otherwise }
\end{array}
\right.
\end{array}
$$
For $i=0$ this follows from Proposition \ref{propertyn} and the fact that 
$\Tor^0_P(\bF(K),\bF)\cong H^*(X_K)$. 
By degree reasons there is no differential ending at or
starting from
$H^{-n}_v(H^0_h(N^*\otimes P)\cong
\Lambda^n(x_1,...,x_n)\otimes \widetilde H^0(K)\cong \widetilde H^0(K)$.
Since $n\geq 2$ this group has total degree $-n+1<0$ and
must therefore
vanish. Hence, the whole column given by $i=0$ and $j\leq -1$ vanishes
and $E_2^{0,0}\cong H^*(X_K)$.
By induction we can apply this argument
successively for $i=1,..., n-2$. In these cases the whole
columns $E_2^{i,j}$ vanish.
For $i=n-1$ there might be a non 
trivial differential $E_2^{n-1,n} \lra  E_2^{0,0}$ and 
we cannot conclude anymore
that $E_2^{n-1, j}=0$.
This shows that 
$\widetilde H^i(K)=0$ for $0\leq i \leq n-2$ and
proves one direction of the claim.

Now we assume that, for $\emptyset\neq \tau\in K$,
all algebras $\bF(\link(\tau))$ are Cohen-Macaulay and that
$\widetilde H^j(K)=0$ for $j < \Dim K=n-1\geq 1$
In particular, $K$ is connected. This implies that
for any pair $\mu, \mu'$ of maximal faces of $K$, there exists a chain
of maximal faces $\mu=\mu_1,...,\mu_s=\mu'$ such that
the intersection $\mu_j\cap \mu_{j+1}$ is non empty. Since $\link{i}$
is pure for all vertices $\{i\}\in K$ (Theorem \ref{proplink}), this implies 
that all maximal
simplices of $K$ have the same order, that $K$ is pure, that
for all faces $\tau \in K$ the dimension of $\link(\tau)$
equals $n-1-\sharp\tau$
and that $H(BT^\tau \times c(\link(\tau)))\cong \bF(\st_K(\tau))$ 
is a finitely generated free 
module over
$P$.  Now we consider again both spectral 
sequences. 
In this case, we get for $H^i_h (H^j_v (N^*\otimes_{P} Q^*))$ 
the same result as above.
And, since $\widetilde H^0(K)=0$, Proposition \ref{propertyn} shows that 
$$
H_v^j(H^i_h(N^*\otimes_{P} Q^*)) \cong 
H_v^j(H^i_h(N^*)\otimes_{P} Q^*)\cong
\left\{
\begin{array}{ll}
\Tor^j_{P}(\bF(K),\bF) & \text{ for } i=0 \\
\Tor^j_{P}(\widetilde H^{n-1}(K),\bF) & \text{ for } i=n \\
0 & \text{ otherwise }
\end{array}
\right.
$$
Hence the
$E_2$-page is concentrated in
two vertical lines given by $i= 0, n$. Considering again total degrees
shows that $\Tor^j_{P}(\bF(K),\bF)=0$  for $j\neq 0$ and that $\bF(K)$
is a finitely generated free $P$-module (Theorem \ref{cmcharfibre}).
This proves the other implication  of the claim.
\qed

\medskip

We can draw the following consequence from the above proof.

\begin{cor} \label{xklimit}
If $\bF(K)$ is Cohen-Macaulay, then there exists a short exact sequence
$$
0 \lra \widetilde H^{n-1}(K) \lra H^*(X_K) \lra \lim H^*(X_{\st(-)}) \lra 
\Lambda^{n-1}(n)\otimes \widetilde H^{n-1}(K)\lra 0 
$$
\end{cor}

\bigskip

\section{Proof of Theorem \ref{gorensteinresult}}

The proof of Theorem \ref{gorensteinresult} is an easy 
consequence of the following statement. 

\begin{thm} \label{gorensteinequivalentstatements}
$\bF(K)$ is Gorenstein${}^*$  if and only if  $\bF(\link_K(\tau))$ is 
Gorenstein${}^*$ 
 for 
all faces $\tau\neq \emptyset$ of $K$ and
$$
\widetilde H^i(K;\bF)\cong
\left\{
\array{ll}
\bF & \text{ for } i=\Dim K \\
0 & \text{ for } i\neq \Dim K
\endarray
\right.
$$
\end{thm}

\medskip

{\it Proof of Theorem \ref{gorensteinresult}:}
We argue as in the proof of Theorem \ref{cmresult}. If $\Dim K=0$, then
$K$ is always Cohen-Macaulay. And $K$ is Gorenstein$^*$ if and only 
if $X_k\simeq S^2$ if and only if $m=2$ if and only if 
$\widetilde H^0(K;\bF)=\bF$.

The general case follows by induction over the dimension of $K$ 
and the above theorem.
\qed

\medskip
 
The proof of Theorem \ref{gorensteinequivalentstatements}
 needs some preparation.

\begin{rem} \label{reduced}
We call a simplicial complex reduced, if for every vertex $i$, the inclusion 
$\st(\{i\}) \subset K$ is proper. If $\st(\{i\}) = K$ 
then $K=\{i\}*\link(\{i\})$.
Hence, any complex $K$ can be written as a
joint product $\Delta*L$ of a full simplex $\Delta$ and a reduced complex $L$.
Moreover, $K$ is reduced if and only if 
the union of all minimal missing faces equals 
the set $V$ of vertices.

We call a simplicial complex $K$ $\bF$-spherical if it satisfies 
the geometric condition of the Gorenstein property, i.e. 
if
$$
\widetilde H^i(\link(\tau);\bF)\cong
\left\{
\array{ll}
\bF & \text{ for } i=\Dim \link(\tau) \\
0 & \text{ for } i\neq \Dim \link(\tau)
\endarray
\right.
$$
for every face $\tau \in K$ in including $\emptyset$.
\end{rem}

\begin{prop} \label{corecondition}
Let $\bF(K)$ be  Gorenstein. Then the following holds:
\nz
(i) $K$ is reduced if and only if $\widetilde H^{n-1}(K)\neq 0$.
In fact, if this is the case, then $\widetilde H^{n-1}(K)\cong \bF$.
\nz
(ii)
If $K$ is reduced, then, for every vertex $i\in V$, 
the link $\link(\{i\})$ is also 
reduced.
\end{prop}

\begin{proof}
By Corollary \ref{xklimit}
we have an exact sequence
$$
0\lra \widetilde H^{n-1}(K) \lra H^*(X_K) 
\lra \prod_{i\in V} H^*(X_{\st(\{i\})})
$$
Let $d\letbe \Dim X_K$. If $\widetilde H^{n-1}(K)=0$ then the second map 
becomes
a monomorphism and there exists an $i\in V$ such that 
$\bF\cong H^{d}(X_K) \lra H^{d}(X_{\st(\{i\})})$ is a monomorphisms.
On the other hand, if  $L\subset K$ is a subcomplex of the same 
dimension as $K$
such that $\bF(K)$ 
and $\bF(L)$ 
are Cohen-Macaulay, then the map 
$H^*(X_K) \lra H^*(X_L)$ is an epimorphism. Hence,
$H^{d}(X_K) \cong H^{d}(X_{\st(\{i\})})\cong \bF$ for some vertex $i\in K$. 
Since both algebras are 
 PD-algebras (Theorem \ref{proplink}), 
this implies that $H^*(X_K) \cong H^{*}(X_{\st(\{i\})})$.
Comparing the two fibrations defining $X_K$ and $X_{\st(\{i\})}$ shows that 
$\bF(\st(\{i\}))\cong \bF(K)$, that $\st(\{i\})=K$, and that $K$ 
is not reduced.

If $K$ is not reduced, then, as the cone of a subcomplex, 
$|K|$ is contractible and 
$\widetilde H^{n-1}(K)=0$. This proves the equivalence in part (i).
If one of the conditions is satisfied, 
then, since $H^*(X_K)$ is $PD$-algebra,
Lemma \ref{topdimensions} shows that 
$\bF\cong H^{n^2+n}(X_K)\cong \widetilde H^{n-1}(K)$.

Since Gorenstein algebras are Cohen-Macaulay, 
Theorem \ref{cmequivalentstatements}
and
Proposition \ref{topdimensions} (ii) tell us that there 
exists a short exact sequence
$$
0\lra H^{n^2+n-2}(X_K) \lra \prod_{i\in V} \widetilde H^{n-2}(\link(\{i\}) 
\lra \widetilde H^{n-1}(K)\cong \bF \lra 0.
$$
 Since $X_K$ is a $PD$-space of dimension $n^2+n$ , the first term 
in the above sequence is isomorphic to $H_2(X_K)\cong \bF^{m-1}$. 
Hence,
the middle term must be isomorphic to $\bF^m$.
Since for each vertex 
$i\in V$, $\Dim_{\bF} \widetilde H^{n-2}(\link(\{i\}))\leq 1$,
this shows that 
$\widetilde H^{n-2}(\link(\{i\}))\cong \bF$ and that $\link(\{i\})$ 
is reduced.
\end{proof}

For $\tau \subset V$ we denote by $K_\tau\subset K$ the full subcomplex 
which consist of all faces $\rho\in K$ 
such that $\rho\subset V\setminus \tau$.
If $\tau=\{i\}$ is a vertex, we denote this complex by $K_i$.

\begin{lem} \label{kicm}
If $K$ is $\bF$-spherical, then, for each vertex $i\in K$, 
the complex $K_i$ is 
Cohen-Macaulay and $\widetilde H^{n-1}(K_i)=0$.
\end{lem}

For the proof we need some preparation.
The inclusions $K_i\subset K$ and $\st(\{i\}) \subset K$ induce epimorphisms
$\bF(K)\lra \bF(K_i)$ and 
$\beta_i:\bF(K)\lra \bF(\st(\{i\}))$ of $P$-modules. 
The kernel of the first map is the ideal $v_i\bF(K)$ generated 
by $v_i\in \bF(K)$. And the second epimorphism induces an isomorphism
$v_i\bF (K) \cong v_i \bF (\st(\{i\}))$. This follows from the 
fact that for any face $\tau\in K$
the monomial $v_i v_\tau=0$ in $\bF(\st(\{i\}))$ 
if and only if $\tau\cup\{i\}\not\in K$.
Moreover, since $\st(\{i\})=\{i\}* \link(\{i\})$, multiplication by $v_i$ 
induces an isomorphism
$\bF(\st(\{i\})) \lra v_i \bF(\st(\{i\})) $. In fact, all the above  maps are 
$\bF(K)$-linear, where $\bF(K)$ acts on $\bF(K_i)$ and $\bF (\st(\{i\}))$ via 
the above projections. Moreover they fit together to a short exact sequence 
$$
0 \lra \bF(\st(\{i\})) \lra \bF(K) \lra \bF(K_i) \lra 0
$$
of $\bF(K)$-modules. The first map is given by $q \mapsto v_i q'$ 
where $\beta_i(q')=q$.
Applying the functor $-\otimes_P$ establishes an epimorphism
$\alpha_i : H^*(X_K) \lra H^*(X_{\st(\{i\})})$ 
and an $H^*(X_K)$-linear map
$$
\psi_i : H^*(X_{\st(\{i\})}) \lra v_i H^*(X_K)
$$ 
given by $\psi_i(a)=v_i a'$ 
where $\alpha_i(a')=a$. In fact, we will show that this map is an isomorphism
(see Corollary \ref{products}).

{\it Proof of Lemma \ref{kicm}:}
Since for two vertices $i,j\in K$ we have $\link_K(\{j\})_i$
equals $\link_K(\{j\})$
if the simplex $\{i,j\}\not\in K$ and equals $\link_{K_i}(\{j\})$
if $\{i,j\}\in K$, we only have to prove that 
$\widetilde H^r(K_i)=0$ for $r\leq n-1$.  
And this claim we prove via an induction over the dimension of $K$.

For $n=1$, $\bF$-spherical implies that $K$ consists
only of two vertices. And for $n=2$,  $\bF$-spherical implies
that $K$ is a triangulation of $S^1$.
In both cases, the claim is straightforward.

Now let us assume that $n\geq 3$.
By excision, 
$H^*(K,K_i)\cong H^*(\st_K(\{i\}),\link_K(\{i\}))\cong 
H^*(\Sigma \link_K(\{i\}))$.
Here $\Sigma \link_K(\{i\})$ denotes the suspension of $\link_K(\{i\})$, 
actually of the geometric realization of $\link_K(\{i\})$.
Moreover, the map $H^{*}(K,K_i) \lra H^{*}(K) $ can be identified with  
$ H^{*}(\Sigma \link_K(\{i\})) \lra H^{*}(K) $ induced by the last 
arrow in the cofibration sequence 
$\link_K(\{i\}) \lra K_i \lra K\lra \Sigma\link_K(\{i\})$.
Since $\Sigma \link_K(\{i\})$ is $\bF$-spherical, it suffices to show that
$H^{n-1}(K_i)=0$.

Let $j\in V$ such that 
$\tau\letbe \{i,j\}\in K$. In the following, 
$K$ will also denote the geometric realization 
of $K$. 
Because of the identities 
$\link_K(\tau)=\link_{\link_K(\{i\})}(\{j\})=\link_{\link_K(\{j\})}(\{i\})$,
$\link_{K_i}(\{j\})=\link_K(\{j\})_i$ and 
$\link_{K_j}(\{i\})=\link_K(\{i\})_j$
all rows and columns in  the homotopy commutative diagram
$$
\CD
\link_K(\tau) @>>> \link_{K_j}(\{i\}) @>>> \link_K(\{i\}) \\
@VVV @VVV @VVV \\
\link_{K_i}(\{j\}) @>>> K_\tau @>>> K_i \\
@VVV @VVV @VVV \\
\link_{K}(\{j\}) @>>> K_j @>>> K
\endCD
$$
consist of  cofibrations.
Passing to suspensions we can extend the diagram 
to the right and to the bottom yielding a homotopy commutative
$4\times 4$-diagram, whose bottom right square looks  like
$$
\CD
K @>>> \Sigma \link_K(\{j\}) \\
@VVV @VVV \\
\Sigma \link_K(\{i\}) @>>> \Sigma^2 \link_K(\tau).
\endCD
$$
In the induced diagram in cohomology in degree $n-1$
$$
\CD
H^{n-1}(\Sigma^2 \link_K(\tau)) @>>> H^{n-1}(\Sigma \link_K(\{j\})) \\
@VVV @VVV \\
H^{n-1}(\Sigma \link_K(\{i\})) @>>> H^{n-1}(K)
\endCD
$$
the left vertical and top horizontal arrows are isomorphisms
by induction hypothesis. Therefore, the other two arrow are either both 
isomomorphisms
or both trivial. And both are isomorphisms if and only if 
$H^{n-1}(K_i)=H^{n-1}(K_j)=0$.

For $n\geq 2$, the complex $K$ is connected and we can connect each pair 
of vertices by 1-dimensional faces. The above argument now shows that
the maps $H^{n-1}(K,K_i)\lra H^{n-1}(K)$ are either 
isomorphisms for all vertices 
or  trivial for all vertices.
Hence it is sufficient to show this map is at 
least nontrivial for at least one vertex, 
or equivalently, that $H^{n-1}(K_i)=0$ for at least one vertex.

Since $H^*(X_K)$ is generated by classes of
 degree 2, a generator of $H^{n^2+n}(X_K) \cong \bF$
can be represented by a monomial 
$a$ which can be written as $v_i a'$
for a suitable vertex $i\in K$. We fix this vertex.
By the above considerations we have 
an exact sequence
$$
H^{n^2+n-2}(X_{\st(\{i\})}) \lra H^{n^2+n}(X_K) \lra H^{n^2+n}(X_{K_i}).
$$
The first map is given by multiplication with $v_i$ 
and therefore an isomorphism.
By Corollary \ref{restrictiontopdegree}, the latter 
map can be identified with
the map $H^{n-1}(K) \lra H^{n-1}(K_i)$, which, 
since $\Dim link_K(\{i\})=n-2$, is an epimorphism. Hence, we have 
$0=H^{n^2+n}(X_{K_i})\cong H^{n-1}(K_i)$, which completes the proof.
\qed

\medskip
 
\begin{cor} \label{products} IF $K$ is $\bF$-spherical, then,
for each vertex $i\in K$, multiplication by $v_i$ induces an isomorphism
$H^*(X_{\st(\{i\})}) \lra v_iH^*(X_K)$.
\end{cor}

\begin{proof}
All terms of the exact sequence 
$$
0\lra \bF(\st(\{i\})) \lra  \bF(K) \lra \bF(K_i) \lra  0
$$
are Cohen-Macaulay (Corollary \ref{proplink}, Theorem \ref{cmresult}, 
Lemma \ref{kicm}. Hence, applying the functor $\otimes_P \bF$ 
establishes a short exact sequences 
$$
0 \lra H^*(X_{\st(\{i\})})\larrow{}  H^*(X_K) \lra H^*(X_{K_i}) \lra 0.
$$
By construction, the first map is given by multiplication by $v_i$.
\end{proof}

\medskip

{\it Proof of Theorem \ref{gorensteinequivalentstatements}:}
If $K$ is Gorenstein${}^*$, then Proposition \ref{corecondition}
shows that $H^{n-1}(K)\cong \bF$ and, together with Corollary \ref{proplink},
that for each face $\tau\in K$ the algebra $\bF(\link(\tau))$ 
is Gorenstein${}^*$.

For the opposite conclusion it suffices to show that
$\soc (H^*(X_K)\cong H^{n^2+n}(X_K)\cong\bF$ 
(Theorem \ref{charfibre}(ii)). By induction we can assume that $K$ is 
$\bF$-spherical.

For each vertex $i\in V$ the map $H^*(X_K) \lra H^*(X_{\st(\{i\})})$ 
is an epimorphism. Hence, this map maps $\soc H^*(X_K)$ to 
$\soc H^*(X_{\st(\{i\})})\cong H^{n^2+n-2}(X_{\st(\{i\})})$. Moreover,
we have an exact sequence 
$0\lra H^{n-1}(K) \lra H^*(X_K) \lra \prod_i H^*(X_{\st(\{i\})})$
(Corollary \ref{xklimit}),
and hence all elements  $\soc H^*(X_K)$ have degree $\geq n^2+n-2$.
Let $a\in H^{n^2+n-2}(X_K)$. This class maps to 
$0\neq b \in H^{n^2+n-2}(X_{\st(\{i\})})\cong \bF$ for some vertex $i\in V$.
And, by Lemma \ref{products}, $0\neq v_i b=v_i a\in H^{n^2+n}(X_K)$. 
This shows that 
$\soc H^*(X_K)\cong H^{n^2+n}(X_K)\cong H^{n-1}(K)\cong \bF$
and finishes the proof.
\qed

%\biskip
%\smallskip

%
%
%
%%
%%
%
%
%
%
%
%
%
%

%%%%%%%%%%%%%%%%%%%%%%%%%%%%%%
%%%%%%%%%%%%%%%%%%%%%%%%%%%%%%
%       Bibliography
%%%%%%%%%%%%%%%%%%%%%%%%%%%%%%
%%%%%%%%%%%%%%%%%%%%%%%%%%%%%%

%\bibliography{algebra,discrete,geometry,topology}

\end{document}